\documentclass[leqno,12pt]{amsart}  
\usepackage{amsmath,amstext,amsthm,amssymb}   
\usepackage[colorlinks=true]{hyperref}
\usepackage{amsrefs}
\numberwithin{equation}{section}

\allowdisplaybreaks
\swapnumbers
\theoremstyle{plain}

\def\bthm#1.#2 #3\ethm{
\begin{\ifx#1ttheorem\fi%
\ifx#1llemma\fi%
\ifx#1ccorollary\fi%
\ifx#1pproposition\fi%
\ifx#1ddefinition\fi}
\label{#1.#2}  
#3 \end{\ifx#1ttheorem\fi%
\ifx#1llemma\fi%
\ifx#1ccorollary\fi%
\ifx#1pproposition\fi%
\ifx#1ddefinition\fi}}

\def\zref #1.#2/{%
\ifx#1e(\ref{e.#2})\fi%
\ifx#1t{Theorem~\ref{t.#2}}\fi%
\ifx#1l{Lemma~\ref{l.#2}}\fi%
\ifx#1p{Proposition~\ref{p.#2}}\fi%
\ifx#1s{Section~\ref{s.#2}}\fi%
\ifx#1c{Corollary~\ref{c.#2}}\fi%
}
\def\Label #1 {\label{#1}}

\def\norm#1.#2.{\lVert#1\rVert_{#2}}
\def\Norm#1.#2.{\bigl\lVert#1\bigr\rVert_{#2}}
\def\NOrm#1.#2.{\biggl\lVert#1\biggr\rVert_{#2}}
\def\NORm#1.#2.{\Bigl\lVert#1\Bigr\rVert_{#2}}
\def\NORM#1.#2.{\Biggl\lVert#1\Biggr\rVert_{#2}}

\def\ip#1,#2.{\langle #1,#2\rangle}
\def\Ip#1,#2.{\bigl\langle#1,#2\bigr\rangle}
\def\IP#1,#2.{\Biggl\langle#1,#2\Biggr\rangle}

\def\abs#1{\lvert#1\rvert}

\newcommand{\za}{\ensuremath{\alpha}}
\newcommand{\zb}{\ensuremath{\beta}}
\newcommand{\zc}{\ensuremath{\psi}}

\newcommand{\ze}{\ensuremath{\epsilon}}
\newcommand{\zve}{\ensuremath{\varepsilon}}
\newcommand{\zvf}{\ensuremath{\varphi}}

\newcommand{\zg}{\ensuremath{\gamma}}
\newcommand{\zG}{\ensuremath{\Gamma}}

\newcommand{\zI}{\ensuremath{\infty}}
\newcommand{\zl}{\ensuremath{\lambda}}
\newcommand{\zk}{\ensuremath{\kappa}}
\newcommand{\zL}{\ensuremath{\Lambda}}
\newcommand{\zm}{\ensuremath{\mu}}

\newcommand{\zw}{\ensuremath{\omega}}
\newcommand{\zW}{\ensuremath{\Omega}}
\newcommand{\zx}{\ensuremath{\xi}}

\newcommand{\zp}{\ensuremath{\pi}}

\def\z#1#2{\ifcase#1 {{\mathcal {#2}}}  
\or {{\mathbf{#2}}}                    
\or { {\boldsymbol{#2}}}                  
\or{{\widetilde{#2}}}                   
\or {{\acute{#2}}}
\or {{\grave{#2}}}
\or {{\bar{#2}}}
\or {\dot{#2}}
\or {\overline{#2}}
\or {{\mathfrak #2}}\fi}
 
\def\ZR{\ensuremath{\mathbb R}}

\def\Xint#1{\mathchoice
   {\XXint\displaystyle\textstyle{#1}}%
   {\XXint\textstyle\scriptstyle{#1}}%
   {\XXint\scriptstyle\scriptscriptstyle{#1}}%
   {\XXint\scriptscriptstyle\scriptscriptstyle{#1}}%
   \!\int}
\def\XXint#1#2#3{{\setbox0=\hbox{$#1{#2#3}{\int}$}
     \vcenter{\hbox{$#2#3$}}\kern-.5\wd0}}

\def\dashint{\Xint-}

 \def\ind#1{{\mathbf 1}_{#1}}

\def\mid{\,:\,}

\def\md#1#2\emd{\ifx0#1
\begin{equation*} #2 \end{equation*}\fi  
\ifx1#1\begin{equation}#2\end{equation}\fi   
\ifx2#1\begin{align*}#2\end{align*}\fi   
\ifx3#1\begin{align}#2\end{align}\fi    
\ifx4#1\begin{gather*}#2\end{gather*}\fi  
\ifx5#1\begin{gather}#2\end{gather}\fi   
\ifx6#1\begin{multline*}#2\end{multline*}\fi  
\ifx7#1\begin{multline}#2\end{mutline}\fi  
}

\def\trans#1{\operatorname{Tr}_{#1}}
\def\modulate#1{\operatorname{Mod}_{#1}}
\def\dilate#1^#2{\operatorname{Dil}_{#1}^{#2}}

\def\sh#1{\operatorname{sh}(#1)}
\def\size #1{\operatorname{size}(#1)}
 \def\Well#1{\operatorname{Well}(#1)}
 \def\emb#1 #2{\operatorname{emb}_{#1}(#2)} 
 \def\Enl#1 #2 {\operatorname{Enl}_{#1} (#2)}

\def\Xint#1{\mathchoice
   {\XXint\displaystyle\textstyle{#1}}%
   {\XXint\textstyle\scriptstyle{#1}}%
   {\XXint\scriptstyle\scriptscriptstyle{#1}}%
   {\XXint\scriptscriptstyle\scriptscriptstyle{#1}}%
   \!\int}
\def\XXint#1#2#3{{\setbox0=\hbox{$#1{#2#3}{\int}$}
     \vcenter{\hbox{$#2#3$}}\kern-.5\wd0}}

\def\dashint{\Xint-}

\def\directions{{\mathsf {drns}}}

\begin{document}

\title[Rubio Littlewood Paley Inequalities]
{Rubio de Francia Littlewood Paley Inequalities 
and Directional Maximal Functions }

\author[G.~Karagulyan]{Grigor Karagulyan }
\thanks{Research supported in 
part by a COBASE grant.} 
\address{Armenian Academy of Sciences
B, Marshall Baghramian Ave.\\ 375019, Yerevan\\ Republic of Armenia
 }

\author[M.T.~Lacey]{Michael T. Lacey}
\thanks{Research supported in part by the 
NSF, and a COBASE grant.}

 \address{School of Mathematics \\ 
  Georgia Institute of Technology \\ 
  Atlanta GA 30332}

  \email{lacey@math.gatech.edu}
 
\date{}

\begin{abstract}  
 In $\ZR^d$, define a maximal function in the directions 
$v\in \directions\subset\{x \mid \abs x=1\}$ by 
\md0 
M^\directions f(x)=\sup_{v\in\directions} \sup_{\zve} \int_{-\ze}^\ze \abs{ f(x-vy)}\; dy.
\emd
For a function $f$ on $\ZR^d$, let $S_\zw f$ denote the Fourier restriction of $f$ to a region
$\zw$. We are especially interested taking \zw\ to be a  sector of $\ZR^d$ with base points at the origin. 
A sector is a  product of the interval $(0,\zI)$
with respect to a choice of   (non orthogonal) basis. What is most important is that the 
basis is a  subset of  $\directions$.  
Consider a collection $\zW$ of pairwise disjoint sectors $\zw$ as above. 
Assume that $M^\directions $ maps $L^p$ into  $L^p$, for some $1<p<\zI  $. 
Then 
we have the following Littlewood--Paley inequality 
\md0
\NORm \Bigl[\sum_{\zw\in\zW}\abs{S_\zw f}^2\Bigr]^{1/2}.q.\lesssim{}\norm f.q.,
\qquad 2\le q<2 \frac p{p-1}.
\emd
The one dimensional analogue of this inequality is due to Rubio de Francia,
\cite{rubio}.  The conclusion when the set of vectors is a fixed  basis is known, is due to 
Journ\'e \cite{MR88d:42028}.  Our method of proof relies on a 
phase plane analysis.  We introduce a notion of Carleson measures 
adapted to $\directions$, and demonstrate a John Nirenberg inequality for these 
measures. The John Nirenberg inequality, and an obvious $L^2$ estimate will prove the Theorem.
\end{abstract}

\maketitle

 \section{Introduction  }  \parskip=11pt

 We are interested in the connection between Littlewood--Paley inequalities in 
 higher dimensions, especially into parallelepipeds and sectors with respect to a variety of distinct 
 bases.  We demonstrate that the  maximal function bounds imply 
 Littlewood--Paley inequalities. 
 
 The classical Littlewood--Paley inequalities concern the decomposition of the frequency
 variables into lacunary pieces.  Our subject is the extent to which these
 inequalities can be generalized when the   
  decompositions of frequency variables  are liberalized.  We specifically
  generalize the beautiful inequality of Rubio de Francia \cite{rubio} to the
  higher dimensions, namely  
  decompositions of frequency variables are specified by an arbitrary collection of pairwise disjoint 
  parallelepipeds and sectors.     The paper concludes with 
  several remarks about our Theorem, its relationship to prior work  and possible generalizations.
  
  In this paper, \zw\ will denote a parallelepiped in $\ZR^d$.  A parallelepiped is
  a product of intervals in  a choice of  (non orthogonal) coordinate axes.  The axes, in particularly, may 
   may vary depending on the parallelepiped.  If the intervals in question are $(0,\zI)$, so that the parallelepiped 
   has a single vertex at the origin, then we say that it is a sector.

  Define the Fourier restriction operator to be 
  \md0
  S_\zw f(x)=\z0F^{-1}\ind\zw \z0F f(x),
  \emd
  where $f$ is a function on the plane and $\z0F f(\zx):=\int_{\ZR^d}e^{-2\zp
  ix\cdot \zx}f(x)\;dx$ is the Fourier transform.  For a collection of parallelepipeds \zW\ 
  set 
  \md0
  S^\zW{}f=\Bigl[\sum_{\zw\in\zW}\abs{S_{\zw}f}^2\Bigr]^{1/2}.
  \emd

  Let $\directions\subset \{x \mid \abs x=1\}$ be a set of norm one vectors in $\ZR^d$. 
  \md0
  M^\directions  f(x)=\sup_{R\in \z0P^\directions } \ind R (x)\dashint_R \abs{f(y)}\; dy
  \emd
  
  \bthm t.lp  Let \zW\ be any collection of disjoint parallelepipeds, and assume that 
  each element of \zW\ is a parallelepiped  with respect to a  basis drawn from 
  vectors in $\directions$.  Assume that $M^\directions $ maps $L^p$ into  $L^p$, for some $1<p<\zI$.
  Then the square
  function $S^\zW$  maps $L^q(\ZR^d)$ into itself for $2\le  q<2 \frac p{p-1}=2p'$.  
  More particularly, for a choice of constant $\zk=\zk(p,d)$, 
  \md1 \label{e.harddirection}
  \norm S^\zW .q.\lesssim{}\norm M^\directions  .p\to p.^{\zk}, \qquad 2\le  q<2 \frac p{p-1}=2p
  \emd
  \ethm
  
  Notice that a sector is an increasing limit of parallelepipeds so that the inequality stated in the abstract is an 
  immediate consequence of the Theorem.

   One would not suppose that the method adopted here would supply an optimal estimate 
for $\zk$ in \zref e.harddirection/.
  
  We adopt a method of proof that emphasizes a notion of space--frequency analysis, 
as following the notes of \cite{laceyrubio}.  That paper concentrates on rectangles  with 
respect to a fixed set of coordinate axes, and the product $BMO$ theory of Chang and R.~Fefferman, 
\cites{MR86g:42038,MR82a:32009}.   More exactly, that paper highlights the role of the 
product Carleson measures in that the case of rectangles with respect to a fixed set of coordinate axes. 

For our current theorem, clearly there is  no such theory,  and so we must find appropriate analogues 
in this setting.

\bigskip 

The one dimensional version of this result is the striking result of Rubio de Francia \cite{rubio}.  
The two dimensional version, with parallelepipeds with respect to a fixed choice of axes, was proved by 
Journ\'e   \cite{MR88d:42028}.  Several other authors have made contributions in this direction, we 
cite, without further comment: Bourgain \cite{MR87m:42008},  C{\'o}rdoba \cite{MR83i:42015}, Olevskii \cites{MR95a:42012a,MR95a:42012b}, Sato \cite{MR92c:42020}, Sj{\"o}lin \cite{MR92c:42020}, and 
Zhu \cite{MR93f:42041}.  These issues are 
surveyed in a recent article by one of us \cite{laceyrubio}.  In particular, the view point we 
take is heavily influenced by that survey article. 

Our theorem, in the case of the plane, and uniformly distributed sectors is due to A.~Cordoba \cite{MR85g:42023}.

\section{Reduction to the Well  Distributed Case} 

This section follows the reduction used by Rubio de Francia \cite{rubio} to collections of intervals 
that are better suited to frequency analysis.   Let $\zW$ be a collection of 
disjoint parallelepipeds \zw, each a parallelepiped with respect to a choice of basis from $\directions$.  Assume in addition that \zW\ satisfies 
\md1\label{e.well} 
\NOrm \sum_{\zw\in\zW} \ind {2\zw}.\zI.<\zI.
\emd
We say that \zW\ is {\em well  distributed.}

\bthm l.well  For any collection $\zW$, we can select a well distributed collection $\Well \zW$ for which we have the inequality 
\md0
\norm S^\zW .q.\lesssim{} \norm M^\directions .(q/2)'.^\zk\norm S^{\Well \zW }.q.,\qquad 1<q<\zI.
\emd
\ethm

In this Lemma, and throughout this paper, $\zk$ denotes a positive number, whose exact value we 
shall not attempt to keep track of.

Consider first the one dimensional case, as we shall be able to pass to higher dimensions by taking 
appropriate products. In turn, in one dimension, we first consider the  interval  $[-\frac12,\frac12]$. Set
\md0
\Well {[-\tfrac12,\tfrac12]}=\{ [-\tfrac1{18},\tfrac1{18}] ,\pm[\tfrac12-\tfrac4{9}(\tfrac45)^{k},\tfrac12-\tfrac4{9}(\tfrac45)^{k+1}]\mid k\ge0\}.
\emd
It is straightforward to check that all the  intervals in this collection have a distance to the boundary of $[-\frac12,\frac12]$ that is 
four times their length. 
In particular, 
this collection is well distributed, and for each $\zw\in \Well {[-\frac12,\frac12]}$ we have 
$2\zw\subset [-\frac12,\frac12]$.    

We define $\Well \zw $ by affine invariance. For an interval \zw, select an affine function $\za\mid [-\frac12,\frac12]\longrightarrow \zw$, we set 
$\Well \zw:=\za(\Well {[-\frac12,\frac12]})$.  And we define $\Well \zW:=\bigcup_{\zw\in\zW}\Well \zw$. It is clear that $\Well \zW $ is well distributed for collections of disjoint intervals $\zW$. 

We shall estimate the the $L^q$ norms of these square functions by duality, as we are always interested in $q>2$.  There is a standard weighted inequality, valid for all $\ze>0$, that  we shall appeal to several times. 
\md1 \label{e.weighted-1}
\int_\ZR  \abs{S_{\zw}f}^2 g\; dx\lesssim\int_\ZR \abs{ S^{\Well {\zw}}f}^2 (M \abs g^{1+\ze})^{1/1+\ze}\; dx ,\qquad   \ze>0.
\emd
Here, $M$ is the one dimensional maximal function.  Clearly, this extends immediately to the 
collection of intervals $\zW$, and the square function $S^\zW$.   

\medskip 

For the case of a parallelepiped $\zw$, we write it as a product of intervals $\zw=\prod_{j=1}^d \zw_j$, in 
the appropriate  non orthogonal coordinates. We then  define $\Well \zw=\prod_{j=1}^d \Well \zw_j$.  By iterating 
\zref e.weighted-1/ in each coordinate, we see that 
\md1 \label{e.weighted-d}
\int_{\ZR^d}  \abs{S_{\zw}f}^2 g\; dx\lesssim\int_{\ZR^d} \abs{ S^{\Well {\zw}}f}^2 (M^\zw \abs g^{1+\ze})^{1/1+\ze}\; dx ,\qquad 1<p<\zI.
\emd
in which $M^\zw$ is a  $d$ times iterate of one dimensional maximal functions in the coordinates associated to \zw.  This Lemma can be summed over a collection of parallelepipeds $\zW$, with the change that the maximal function $M^\directions $, iterated $d$ times,  must be imposed on the right hand side of the inequality.   Thus  the Lemma is proved. 

\medskip  

We use the notion of well distributed just as Rubio de Francia did, to pass from the sharp Fourier restriction, with it's accompanied long range spatial behavior, to convolution with Schwartz functions, 
with very rapid spatial decay.  

\bthm  l.smooth    We assume that $\zW$ is a set of parallelepipeds in the set of directions $\directions$.  
Assume that $M^\directions$ is bounded on $L^p$.  Let $\zc_\zw$, for $\zw\in\zW$, be a Schwartz function with $\ind {\zw}\le\widehat{\zc_\zw}\le \ind {2\zw}$.  Then, 
\md0
\norm S^\zW f .q.\lesssim\norm M^\directions .p\to p.^\zk \NOrm \Bigl[ \sum_{\zw\in\zW} 
\abs{\zc_\zw*f }^2 \Bigr]^{1/2} .q.,\qquad 2<q<2p'.
\emd
\ethm

The proof is again an  application of \zref e.weighted-d/, in exactly the same manner.

   \section{The Space--Frequency Tiles}
 Our purpose is to define a discrete analog of the square function $S^\zW$.  The discrete analog will more easily permit an analysis in terms of the theory of  Carleson measures in directions  that we develop in 
 \zref s.carleson/ and the subsequent sections.

Let us introduce the operators associated to translation, modulation and dilation. 
\md5\label{e.translation} 
\trans y f(x)=f(x-y), \quad y\in\ZR^d
\\  \label{e.modulate}
\modulate \xi f(x)=e^{i \xi \cdot x}f(x),\quad \zx\in\ZR^d.
\emd

A parallelepiped $R\subset\ZR^d$ is a product of intervals $R=\prod_{i=1}^d R_i$ with respect to a choice of  basis 
in $\ZR^d$, with coordinates $(x_1,x_2,\ldots,x_d)$. Set a dilation operator associated to $R$ to be 
\md1 \label{e.dilate}
\dilate R^p f(x)= \abs{R}^{-1/p}f(\za(x)),\qquad 0\le{}p<\infty, 
\emd  
Here, $\za$ is an affine map that carries $R=\prod_{i=1}^d R_i$ coordinate wise into the standard cube $ [-1/2,1/2]^d$.  
This depends on the orientation of $R$, and its side lengths of $R$,  and its  location.    Notice that this definition implictly 
incorporates a translation. 

Two parallelepipeds $\zw$ and $R$ are said to be {\em dual} if they are both parallelepipeds with respect to the same choice of coordinate 
vectors, and with those vectors ordered, writing $\zw=\prod_{j=1}^d\zw_j$ and $R=\prod_{j=1}^d R_j$, one has 
\md1\label{e.dual} 
1\le{}\abs{R_j}\cdot\abs{\zw_j}\le2,\qquad 1\le{}j\le{}d.
\emd

Let $\z0D(\directions)$ denote the collection of parallelepipeds that are dyadic with respect to a choice of basis  in $\directions$.  That is, the parallelepiped is a product of dyadic intervals in it's basis.\footnote{This restriction is made to keep $\z0D(\directions)$ countable.  Our Carleson measure theory can then be phrased in terms of sums, rather than integrals. }

Call a product $R\times \zw$ a {\em tile} if $\zw$ and $R$ are dual, and $R\in\z0D(\directions)$.
We shall associate to each tile an appropriate function.   
Fix a Schwartz function $\zvf\ge0$ with $\widehat\zvf$ supported on $[-9/16,9/16]^d$, set 
\md1 \label{e.f_s.def}
\zvf_{R\times \zw}={}  \modulate {c(\zw)} \dilate R^2 \zvf
\emd
 where $c(J)$ denotes the center of $J$.  

Let $\z0T(\zW)$ be a collection of tiles so that for all $s=R_s\times \zw_s\in\z0T(\zW)$,  we have that $\zw_s\in\zW$. 
Our space--frequency square function is 
\md1 \label{e.sf}
SF^\zW{}f=\Bigl[\sum_{s\in\z0T(\zW)} \frac{\abs{\ip f,\zvf_s.}^2}{\abs{R_s}}\ind {R_s} \Bigr]^{1/2}.
\emd
We let $\directions$ be the set of coordinates for the parallelepipeds in \zW.  

\bthm p.well  Assume that $M^\directions $ maps $L^p$ into  $L^p$ for some 
$1<p<\zI$.  Then,  for all well distributed $\zW$, we have 
\md0
\norm SF^\zW.q.\lesssim \norm M^\directions.p.^\zk ,\qquad 2\le{}q<2p'.
\emd
\ethm

We impose the well distributed assumption to trivialize the boundedness of the square function on $L^2$.  
Indeed, by our construction, the function $\zvf_s$ is supported on $2\zw_s$, and these sets have bounded overlap.  Thus, the $L^2$ inequality reduces to checking it for $\zW$ consisting of just one parallelepiped.  But then, the fact that $\zvf$ is a Schwartz function proves the desired inequality.

 The main point in the proof of the proposition is this fact, which we state in the language of 
 \zref s.carleson/.

 \bthm l.cm  Fix $\ze>0$. Assume that $M^\directions $ maps $L^p$ into  $L^p$ for some 
$1<p<\zI$, and that $\zW$ is well distributed.  Then for each $f$ on $\ZR^d$, with $L^\zI$ norm bounded by one, the map from $\z0D(\directions)$ to $\ZR_+$ below 
has $CM(\directions)$  norm at most $\norm M^\directions .p\to p.^{\zk}$. 
\md0
\z0D(\directions)\ni R\longrightarrow \sum_{\substack { s\in\z0T(\zW)\\ R_s=R }} \abs{ \ip f,\zvf_s .}^2 
\emd
\ethm

 \section{Proofs}

\begin{proof}[Proof of \zref p.well/]
Let $q=2\frac p{p-1}=2p'$.  As the $L^2$ inequality is a consequence of the well distributed assumption, we need only show the restricted weak type inequality at $L^{q}$. 
Fix a function $f=\ind F$, where $F$ is of finite measure.  
For a set of of tiles $\z0S\subset \z0T(\zW)$, let 
\md0
S(\z0S):=\sum_{s\in\z0S} \frac{\abs{\ip f,\zvf_s.}^2}{\abs{R_s}} \ind {R_s} 
\emd
[Note the lack of a square root here.] It is our task to show that for all $\zl>0$, 
\md1 \label{e.2222} 
\abs{ \{ S(\z0T(\zW))>\zl\}}\lesssim{}\zl^{-p'}\norm M^\directions .p.^{\zk}  \abs F .
\emd
As the admissible collections of tiles are invariant under  dilations which are uniform 
in all coordinates, it suffices to consider the case of $\zl=1$ in this inequality. 

In addition define 
\md4
\sh {\z0S} {}:={} \bigcup_{s\in\z0S}R_s.
\\
\size {\z0S} {}:={}\sup_{\z0S'\subset\z0S} \Bigl[ \abs{ \sh {\z0S' } }^{-1}\sum_{s\in\z0S'} \abs{\ip f,\zvf_s.}^2\Bigr]^{1/2}
\emd
The definition of size is equivalent to that of a Carleson measure. 
It follows immediately from \zref l.cm/ that we have $\size {\z0T(\zW)}   {}\le{}
\zm\lesssim\norm M^\directions .p.^\zk$.   And moreover, from                                                                                          
 \zref l.jn/   we have the inequality of John Nirenberg type 
\md1  \label{e.Sjn} 
\norm S(\z0S).p'.\lesssim{}   \size {\z0S}^2
\abs{\sh{\z0S}}^{1/p'}.
\emd

The next stage of the argument is the primary decomposition of the set 
of tiles $\z0T(\zW)$.  
Note that if there are two disjoint subcollections $\z0S^j$, $j=1,2$ of $\z0T(\zW)$ such 
that  for both we have 
\md0
\sum_{s\in\z0S^j}\abs{\ip \ind F, \zvf_s. }^2\ge\frac12 \abs{\sh {\z0S^j}}\zm^2,
\emd
then  the same inequality holds for their union. Thus there is a maximal (not necessarily unique) subcollection $\z0S_1\subset \z0T(\zW)$ such that 
\md4
\sum_{s\in\z0S_1}\abs{\ip \ind F, \zvf_s. }^2\ge\tfrac12 \abs{\sh {\z0S_1}}\zm^2,
\\
\size {\z0T(\zW)-\z0S_1}\le  \abs{\sh {\z0S}}\zm.
\emd
The last part suggests a recursive application to $\z0T(\zW)-\z0S_1$.  Carrying this out will result in 
a decomposition of $\z0T(\zW)$ into collections $\z0S_k$, $k\ge-2\log \zm$, for which for each $k$, we have 
\md5  \label{e.sk>} 
\sum_{s\in\z0S_k}\abs{\ip \ind F, \zvf_s. }^2\ge2^{-2k}\zm^2\abs{\sh {\z0S_k}},
\\  \nonumber
\size {\z0S_k}\lesssim2^{-k}\zm	.
\emd
The top inequality gives an upper bound on $\abs{\sh {\z0S_k}}\lesssim{} 2^{2k}\zm^2\abs F $. This follows from  the $L^2$ inequality that holds by design. 

Then, for a small value of $\ze$, we will have $\sum_{k\ge -2\log \zm} \ze 2^{-\ze k}\le1$, so that  it suffices to see that 
\md0
\sum_{k\ge -2\log \zm} \abs{\{ S(\z0S_k)>\ze 2^{-\ze k}\}}\lesssim{}\zm^\zk \abs F ,
\emd
for some absolute choice of $\zk$. 
We can employ the  John Nirenberg inequality \zref e.Sjn/ to see that 
\md2
\abs{\{ S(\z0S_k)>\ze 2^{-\ze k}\}}\le{}& 2^{-k(-\ze)2p'}\abs {\sh {\z0S_k}}
\\{}\lesssim{}& 2^{-k[(1-\ze)2p'-2]}\abs F .
\emd
This is summable in $k\ge -2\log \zm$ to $\zm^\zk \abs F $.  
\end{proof}


\begin{proof}[Proof of \zref l.cm/.]  
The main points are the well distributed assumption, and  the boundedness of the maximal function $M^\directions$.
Fix a function $f$ on $\ZR^d$ that is bounded by one, and  an open set $U\subset\ZR^d$. 
We are to show that 
\md0
\sum_{\substack{ s\in\z0T(\zW) \\ R_s\subset U}}\abs{\ip f,\zvf_s.}^2\lesssim{} \norm M^\directions .p\to p.^\zk \norm f.\zI.^2\abs U .
\emd
We need this elementary Lemma.

  \bthm l.elem   For all $N\ge1$, $0<a<1$,   $\zw\in\zW$ and functions $f$ supported off of the set $\{ M^\directions \ind U >a^d\}$, we have 
  the estimate 
  \md0
  \sum_{\substack{R_s\subset U \\\zw_s=\zw}}\abs{\ip f,\zvf_s. }^2\lesssim{}a^N \norm f.2.^2
  \emd
  \ethm 
  
  \begin{proof} 
  Let $\z0T$ consist of those tiles $s\in\z0T(\zW)$ for which $R_s\subset U$ and $\zw_s=\zw$.   Let 
  $\zG_s=\frac15a^{-1}R_s$, and note that this parallelepiped cannot intersect the support of $f$.    We decompose 
  $\zvf_s=\za_s+\zb_s$, where $\za_s$ is a smooth function supported on $\zG_s$, and equal to $\zvf_s$ on $\frac12 \zG_s$. 
  
  Thus, $\zb_s$ is the ``trivial" part.  And indeed, it is straightforward to verify that 
  \md0
  \sum_{s\in\z0T} \abs{\ip g,\zb_s.}^2\lesssim{} a^N  \norm g.2.^2.
  \emd
  But we will not apply this estimate to $f$, but rather to $\zc_{2\zw}*f$.  Then, we have 
  \md0
  \ip f,\zvf_s.=\ip \zvf^{2\zw}*f,\zvf_s.=\ip  \zvf^{2\zw}*f,\za_s.+ \ip \zvf^{2\zw}*f,\zb_s.
  \emd
  
  Certainly, we do not need to futher consider the inner products with $\zb_s$.   As for the inner products with $\za_s$, the 
  main point is that $ \zvf^{2\zw}*f$ is dominated by appropriate iterate of a  maximal function $M$ in coordinates for the parallelepiped $\zw$.  And in fact, using the assumption about the 
  support of $f$, and the fact that $\zvf$ is a Schwartz function, for each $s\in\z0T$, and $x\in \zG_s$, we have 
  \md0
  \abs{ \zvf^{2\zw}*f(x)}\lesssim{} a^N M f(x). 
  \emd
  Therefore, we can estimate by Cauchy--Schwarz, 
  \md2 
  \sum_{s\in\z0T} \abs{\ip  \zvf^{2\zw}*f,\za_s.}^2\lesssim{}&a^N \sum_{s\in\z0T} \int_{\zG_s} \abs{ M f }^2\; dx 
  \\{}\lesssim{}& a^{N-d} \norm f.2.^2.
  \emd
  The last line follows from the  $L^2$ bound on $M$ that is uniform in the choice of $\zw$, and  since the rectangles $\zG_s$ can overlap at most ${}\lesssim a^{-d}$ times. 
  
  \end{proof}

  Let us decompose $f$ into a sum of functions $f_k$, for $k\ge0$, 
  \md4 
  f_0=f\ind { \{ 1/2\le{} M^\directions \ind U \} } , 
  \\
  f_k=f\ind { \{ 2^{-k-1}\le M^\directions \ind U \le 2^{-k}\}}, \qquad k>0.
  \emd
  For $f_0$, we use the well distributed assumption, and the Bessel inequality it implies, to see that 
  \md0
  \sum_{
\substack{s\in\z0T(\zW) \\R_s\subset U}} \abs{\ip f_0,\zvf_s.}^2\lesssim\norm f_0.2.^2\le\norm f.\zI.^2\abs U   
  \emd
  
  The terms $f_k$ require \zref  l.elem/.  For each $\zw\in\zW$, 
  \md0
  \sum_{\substack{ R_s\subset U  \\ \zw_s=\zw}} \abs{\ip f_k,\zvf_s.}^2\lesssim{} 2^{-kp}\norm \zc_\zw *f_k.2.^2
 \emd 
 By the well distributed assumption, this estimate can be summed over $\zw$, to acheive the bound 
 \md0
 \sum_{R_s\subset U}\abs{\ip f_k,\zvf_s.}^2\lesssim{} 2^{-kp}\norm  f_k.2.^2
 \emd
 The last fact to be noted is that by the boundedness of the strong maximal function 
 on $L^p$, we have that $\norm f_k.2.\lesssim{} 2^{kp/2}\norm f.\zI.$.   Therefore, this estimate can be summed over 
 $k>0$, to conclude the proof. 
 \end{proof}

\begin{proof}[Proof of \zref t.lp/.]  We indicate the proof of  our main Theorem.  Recall from \zref l.well/ that it suffices to consdier well distributed collections $\zW$, and so we  should argue from \zref p.well/, and 
construct a smooth square function, as in \zref l.smooth/, for which we have the same norm inequalities 
as in \zref p.well/.   Namely, if the maximal function $M^\directions$ is bounded on $L^p$, we should 
obtain norm inequalities for the square function for $2<q<2p'$. 
The argument that we present here is comprised of standard lines of reasoning. 

Let us set 
\md4
B^{\zW} f:=\Bigl[ \sum_{\zw\in\zW} \abs{B_\zw f}^2  \Bigr]^{1/2},
\\
B_\zw{}f:=\sum_{\substack{s\in\z0T(\zW)\\ \zw_s=\zw}} \ip f, \zvf_s. \zvf_s 
\emd
We observe that $B^\zW$ maps $L^q$ into itself for $2<q<2p'$.    Indeed, one has the inequality 
\md2 
B^{\zW}f\lesssim{} & \Bigl[   \sum_{{s\in\z0T(\zW)}}  \abs{ \ip f,\zvf_s. \zvf_s }^2 \Bigr]^{1/2}
\\{}\lesssim{}& \Bigl[   \sum_{{s\in\z0T(\zW)}}  \frac{\abs{ \ip f,\zvf_s. \zvf_s }^2}{\abs{R_s} }
(M^{R_s} \ind {R_s})^2\Bigr]^{1/2}.
\emd
Here, in the top line, we are using Cauchy---Schwarz and the rapid decay of the functions $\zvf_s$.  
Notice that the bottom line is very similar to the square function considered in \zref p.well/.  The only 
difference is the imposition of the maximal function on the indicator of $R_s$.   In particular, $M^{R_s}$ is the maximal function in the basis for the rectangle $R_s$. 

There is a counterpart of \zref e.weighted-1/ that applies to the maximal function.  Namely, in one dimension, 
\md0
\int \abs{M f}^2g\; dx\lesssim{} \int \abs f (M\abs g^{1+\ze})^{1/1+\ze}\; dx,\qquad  \ze>0.
\emd
Applying this to the right hand side above, and using the hypothesis that $\norm M^\directions.p\to p.$ is finite, we see that $B^\zW$ satisfies the claimed range of $L^q$ inequalities.  

\medskip 

From $B^\zW$ we should pass to an operator which is a square function  of convolution operators as in 
\zref l.smooth/.  It suffices to consider each parallelepiped $\zw\in\zW$ individually. Consider the limit 
\md0
C_\zw f:=\lim_{Y\to\zI} \int \trans {-y} B_\zw \trans y \; \zm_Y(dy)
\emd
where $\zm_Y$ is normalized Lebesgue measure on the ball centered at the origin of radius $Y$.  
One sees that this limit exists for all Schwartz functions.  Since $B_\zw$ is clearly a bounded operator on $L^2$, we conclude that $C_\zw$ is as well.  It also commutes with all translation operators, by construction. Hence, it is a convolution operator.  And one may check directly that $C_\zw f=\zc_\zw$, where 
\md0
\zc_\zw(x)=\int_{\ZR^d}\zvf(x+y)\overline{\zvf_\zw(y)}\; dy
\emd
By construction, $\ind \zw\le{}\widehat{\zc_\zw}\le{}\ind {2\zw}$, so that we have constructed a smooth square function as in \zref l.smooth/.  The proof of the Theorem is complete. 
\end{proof}


\section{Carleson Measures with Directions}  \label{s.carleson}

We set out a theory of Carleson measures associated to sets of directions $\directions$ in $\ZR^d$.  
Recall that $\z0D(\directions)$ denotes the set of parallelepipeds in $\ZR^d$ that are dyadic with respect 
to a choice of bases from $\directions$.   When $\directions$ is a single  orthogonal basis, all of this reduces to the 
Carleson measure theory associated with product $BMO$.

For a function $\zL\mid \z0D(\directions)\longrightarrow \ZR_+$, we 
set 
\md1 \label {e.defvcm}
\norm \zL. CM(\directions).:=\sup_{U\subset \ZR^d} \abs U^{-1}\sum_{R\subset U}\zL(R)
\emd
What is to be emphasized, is that the supremum is taken over all subsets of $\ZR^d$ of 
finite measure.

Despite the generality of these definitions, it does permit the development of a rudimentary theory. 
The first fact to note is an extension of the John Nirenberg inequality. 

\bthm l.jn  Assume that $M^\directions $ maps $L^p$ into weak $L^p$.  Then we have the inequality 
below, valid for all sets $U\subset \ZR^d$, of finite measure. 
\md0
\NOrm \sum_{R\subset U}\frac{\zL(R)}{\abs R} \ind R .q.\lesssim \norm M^\directions  .p\to p,\zI. \norm \zL.CM(\directions). ,\qquad 1\le{}q\le{} p'. 
\emd
\ethm 

\begin{proof}
The argument of \cite{MR82a:32009} for the John Nirenberg inequality needs only modest 
modifications in the present setting.  
 Define 
 \md0
 F_U=\sum_{R\subset U}\frac{\zL(R)}{\abs R} \ind R
 \emd
 We want to show that $\norm F_U.p'.\lesssim{}\norm M^\directions  .p\to p,\zI.\abs{U}^{1/p'}$.  This we shall do by
 showing that there is an open  $V\subset\ZR^d$ so that $\abs V<\frac12\abs U$ so that 
 \md0
 \norm F_U .p'.\lesssim{}\norm M^\directions  .p\to p,\zI.\abs{U}^{1/p'}+\norm F_V.p'.
 \emd
 An inductive argument proves the desired inequality. 
 
 This is done by duality.  Thus, choose $g\in L^p$ of norm one so that
  $ \norm F_U .p'.=\ip
 F_U,g.$.  Then let $V=\{M^\directions  g>c\abs{U}^{-1/p}\}$.  For appropriate constant $c\simeq{}\norm M^\directions  .p\to p,\zI.$, the measure of $V$ is at most half of the measure of $U$.  Note that 
 if $R\not\subset V$, then $\dashint_{R}g\;dx<c\abs{U}^{-1/p}$. Hence, 
  \md2
 \norm F_U .p'.={}&\ip F_U,g.
 \\
 {}={}&\sum_{\substack{R\subset U\\R\not\subset V}}\zL(R) \dashint_R g \; dx+{} 
 \ip F_V, g.
 \\{}\lesssim{}&\norm M^\directions  .p\to p,\zI. \abs{U}^{1/p'}+\norm F_V.p'.
 \emd

 \end{proof}

 
 \section{Concluding Remarks} 
 
 Our result is unsatisfying, as it does not give new examples of maximal functions $M^\directions$ for which there are $L^p$ bounds.   Indeed, the only instance in which we have a range of examples is the 
 plane, and it is this setting that these results will likely find application.

 It would also be of interest to know that the theorem is sharp as to the range of indicies that we prove the 
 square function inequalities.  We can show that it is sharp in one case, in that of uniformly distributed directions in the plane. 
 
For a large integer $N$, let $\zG$ be a collection of pairwise disjoint sectors in the plane, with vertexes at the origin, and 
opening angle $2\zp/N$.    The maximal function we have associated to this square function is one over two sets of uniformly distributed directions in the plane.  It is well known that this maximal function 
admits  bounds that are logarithmic in $N$ for $p\ge2$, but the bound blows up as a power of $N$ for $1<p\le2$.  We conclude that the square function $S^\zG$ admits logarithmic bounds for $2<q<4$.  

A simple example shows that this is the correct range of indicies for which one has such a bound.  Let 
$\zvf$ be a Schwartz function with $\widehat\zvf$ non negative, radial, rotationally symmetric, and supported in a small 
annulus about $\abs \zx=1$.  For each $\zg\in\zG$, it is routine to see that $\abs{S_\zg \zvf}\gtrsim{} 
N^{-1} \ind {R_\zg } $, where $R_\zg$ is a $1\times N$ rectangle, with center at the origin, and 
long direction oriented in the direction of the bisectrix of $\zg$.  

We therefore have 
\md2
\norm S^\zG \zvf .q.^2\gtrsim{}& \Norm N^{-1}\sum_{\zg\in\zG} \ind {R_\zg}  .q/2.
\\
{}\gtrsim{}&N^{1-4/q},\qquad 4<q<\zI,
\emd
preventing the possibility of a meaningful result for the square function on $L^q$ for $q>4$.

\end{document}